\pdfoutput=1
\documentclass[11pt,a4paper]{article}
\usepackage[T1]{fontenc}
\usepackage[utf8]{inputenc}
\usepackage{lmodern}
\usepackage[margin=2.7cm]{geometry}
\usepackage{amsmath,amssymb,amsthm,mathtools}
\usepackage{booktabs,microtype}
\usepackage{enumitem}
\usepackage{tikz}
\usepackage{pgfplots}
\pgfplotsset{compat=1.18}
\usepackage[authoryear,round]{natbib}
\usepackage[colorlinks=true,linkcolor=red,citecolor=blue,urlcolor=blue]{hyperref}
\newtheorem{theorem}{Theorem}[section]
\newtheorem{proposition}[theorem]{Proposition}
\newtheorem{corollary}[theorem]{Corollary}

\theoremstyle{definition}
\newtheorem{definition}[theorem]{Definition}
\newtheorem{remark}[theorem]{Remark}
\newtheorem{objective}[theorem]{Objective}
\newcommand{\R}{\mathbb{R}}
\newcommand{\eps}{\varepsilon}
\newcommand{\calE}{\mathcal{E}}
\newcommand{\de}{\Delta}
\title{A Feasible-Domain Reduction Approach with Objective-Loss Guarantees for Multilevel Multiobjective Linear Programming}
\author{Mustapha Kaci\\[0.4em]
{\small Laboratory of Mathematics and Statistics and Their Applications, USTOMB.}\\
{\small Djillali Liabes University of Sidi Bel Abbes, Sidi Bel Abbes, Algeria.}\\[0.3em]
{\small \href{mailto:kaci.mustapha.95@gmail.com}{\nolinkurl{kaci.mustapha.95@gmail.com}}\quad and\quad
\href{mailto:mustapha.kaci@univ-sba.dz}{\nolinkurl{mustapha.kaci@univ-sba.dz}}}}
\hypersetup{pdfauthor={Mustapha Kaci}}
\date{}
\begin{document}
\maketitle
\begin{abstract} In some approaches to decentralized multilevel multiobjective programming, each decision-maker selects a compromise and restricts the feasible region before passing the decision to the next level. However, these restrictions do not necessarily limit the deterioration of objective values at preceding levels. This paper studies how to choose variable bounds that keep such losses within prescribed tolerances. We establish exact conditions using linear programming duality and derive simpler bounds based on objective coefficients. We also examine the choice of restriction widths and establish guarantees for successive reductions. A two-level numerical example illustrates the results and the influence of tolerances on the final compromise. \end{abstract}


\section{Introduction}\label{sec:intro}

In multilevel multiobjective linear programming, several decision-makers
act in a prescribed order. Each chooses some of the variables and evaluates
the decision through several criteria. A restriction imposed at one level
changes the choices available to the next. It can also reduce the values
of objectives that earlier decision-makers sought to preserve.

Different approaches to hierarchical decisions have been proposed.
\citet{AnandalingamApprey1991} studied conflict resolution, and
\citet{ShihLaiLee1996} developed a fuzzy approach. Successive restrictions
of decision variables appear in \citet{SinhaSinha2004}; cooperative and
noncooperative formulations are discussed by \citet{SakawaNishizaki2009}.
Related bilevel models, including those with vector objectives, are studied
in \citet{DempeZemkoho2020,DempeMehlitz2025,LafhimZemkoho2023,HoffMehlitz2024}.

The geometry of efficient solutions in linear multiobjective optimization
has been studied by \citet{YuZeleny1975} and \citet{Benson1998}; see also
\citet{Ehrgott2005}. Prescribed objective thresholds are familiar from the
$\eps$-constraint method \citep{Mavrotas2009}. Here they serve a different
purpose: we use them to check whether \emph{every decision} left in a reduced
region respects an earlier decision-maker's tolerances.

\citet{KaciRadjef2022} described common compromises through nondominated
faces of the feasible polytope. Their adaptive method
\citep{KaciRadjef2023} uses successive bounds on decision variables.
However, the widths of these bounds do not state how much earlier objectives
may deteriorate. A preliminary study of this question appears in
\citet{Kaci2026Note}.

We give exact tests for the largest possible loss on a restricted region and
simpler tests using only objective coefficients. We then determine admissible
restriction widths, prove that previous guarantees persist after further
reductions, and derive perturbation bounds. The results are illustrated
on a two-level sorting segment.

We assume that a suitable sorting polytope has already been found. Our aim
is to control compromise selection within that polytope, not to characterize
solutions of the original nested best-response problem.

Section~\ref{sec:prelim} recalls the notation from
\citet{KaciRadjef2022}. Section~\ref{sec:model} states the problem.
Sections~\ref{sec:cert} and~\ref{sec:parametric} give loss tests and the
choice of bounds. Section~\ref{sec:hierarchy} considers successive reductions,
Section~\ref{sec:stability} establishes stability estimates, and
Section~\ref{sec:numerical} presents the example.
Section~\ref{sec:conclusion} concludes the paper.

\section{Preliminaries}\label{sec:prelim}

We use the notation of the geometric approach in \citet{KaciRadjef2022}.
There are $P\geq2$ decision-makers. The variables are divided into blocks,
\begin{equation}\label{eq:partition}
x=(x^{(1)},\ldots,x^{(P)})\in\R^n,\qquad
x^{(p)}\in\R^{n_p},\qquad \sum_{p=1}^P n_p=n,
\end{equation}
with $x^{(p)}$ controlled by decision-maker $p$. The common feasible set is
the nonempty compact polytope
\begin{equation}\label{eq:S}
S=\{x\in\R^n:Ax\leq b,\ x\geq0\},
\end{equation}
where $A\in\R^{m\times n}$ and $b\in\R^m$.

At level $p$, the $k_p\geq1$ criteria are the components of
\begin{equation}\label{eq:F}
F_p(x)=C_px+d_p\in\R^{k_p},\qquad
C_p=\begin{bmatrix}
c_{p1}^{\mathsf T}\\[-2pt]\vdots\\c_{pk_p}^{\mathsf T}
\end{bmatrix}.
\end{equation}
\subsection{Efficient solutions and common compromises}

For $u,v\in\R^{k_p}$, the inequality
$u\geq v$ is understood componentwise.
A feasible point $z$ is efficient for
decision-maker $p$ if
\begin{equation}\label{eq:eff}
\nexists\,y\in S:
F_p(y)\geq F_p(z),\qquad
F_p(y)\neq F_p(z).
\end{equation}

Let $E_p$ denote the set of efficient
solutions of $F_p$ over $S$.

Following \citet{KaciRadjef2023},
the set of possible common compromises is
\begin{equation}\label{eq:E}
\calE=\bigcap_{p=1}^{P}E_p.
\end{equation}

Thus, every point of $\calE$ is efficient
for all decision-makers. This set is not
necessarily convex.

\subsection{Sorting polytope}

As shown in \citet{KaciRadjef2022},
the common efficient set can be described
through non-dominated faces of $S$.
These faces provide convex regions for
the selection of a compromise.

We consider a nonempty compact convex
polytope $T$ such that
\begin{equation}\label{eq:T}
\varnothing\neq T\subseteq\calE.
\end{equation}
We call $T$ a sorting polytope and
assume that it has already been determined.

In accordance with the face-based
construction, we assume that $T$ is
a face of $S$. It can then be written as
\begin{equation}\label{eq:Trep}
T=\left\{x\in S:
A_{i\cdot}x=b_i\ (i\in I),\quad
x_j=0\ (j\in J)\right\},
\end{equation}
where $I\subseteq\{1,\ldots,m\}$ and
$J\subseteq\{1,\ldots,n\}$.

For subsequent linear programs, define
\begin{equation}\label{eq:Tineq}
\widehat A=
\begin{bmatrix}
A\\
-I_n\\
-A_{I\cdot}\\
Q_J
\end{bmatrix},
\qquad
\widehat b=
\begin{bmatrix}
b\\
0_n\\
-b_I\\
0_{|J|}
\end{bmatrix},
\end{equation}
where $Q_J$ consists of the rows
$e_j^{\mathsf T}$ for $j\in J$.

Hence,
\[
T=\{x\in\R^n:
\widehat Ax\leq\widehat b\}.
\]

This representation involves only the
inequalities defining $S$ and the
equalities corresponding to the selected
active constraints.

\section{Problem formulation}\label{sec:model}

With the preceding notation, the hierarchical model is
\begin{equation}\label{eq:MLMOLP}
\left\{
\begin{aligned}
\text{Level }1:\quad
&\max_{x^{(1)}} F_1(x)=C_1x+d_1,\\
&\text{where }(x^{(2)},\ldots,x^{(P)})\text{ solves}\\[2mm]
\text{Level }2:\quad
&\max_{x^{(2)}} F_2(x)=C_2x+d_2,\\
&\text{where }(x^{(3)},\ldots,x^{(P)})\text{ solves}\\
&\vdots\\
\text{Level }P:\quad
&\max_{x^{(P)}} F_P(x)=C_Px+d_P,\\
&\text{subject to }x\in S.
\end{aligned}
\right.
\end{equation}

\begin{objective}[Objective of the study]
    In the adaptive method of \citet{KaciRadjef2023}, a decision-maker selects
a compromise and then restricts the variables available to the following
levels. Such restrictions keep the decision within the chosen region, but
do not give a bound on the loss in earlier objectives. We seek bounds that
do provide this guarantee.
\end{objective}

Let $z\in T$ and define
\begin{equation}\label{eq:D}
D=\{x\in T:\ell\leq x\leq u\}.
\end{equation}

\begin{definition}[Objective-retaining restriction]
Let $\eps_p\in\R_+^{k_p}$ be the prescribed
tolerances of decision-maker $p$.
The restriction $D$ is said to preserve
the objectives of decision-maker $p$
within $\eps_p$ if
\begin{equation}\label{eq:retention}
F_p(x)\geq F_p(z)-\eps_p,
\qquad \forall x\in D.
\end{equation}
\end{definition}

Since $F_p$ is affine, this condition is
equivalent to
\begin{equation}\label{eq:boxguard}
C_px\geq C_pz-\eps_p,
\quad
\forall x\in D.
\end{equation}

Thus, every decision remaining in $D$
satisfies the prescribed objective tolerances.

We next give tests that decide whether this condition holds, before the
reduced domain is passed to the next decision-maker.

\section{Main results}
\subsection{Checking objective losses}\label{sec:cert}

We first measure the largest deterioration of one criterion over a proposed
domain.  Put
\begin{equation}\label{eq:Delta}
\de_{pj}(D;z)=\max_{x\in D}c_{pj}^{\mathsf T}(z-x)
=c_{pj}^{\mathsf T}z-\min_{x\in D}c_{pj}^{\mathsf T}x.
\end{equation}

\begin{theorem}[Exact loss test]\label{thm:exact}
Let $z\in D$. Then, for each $p$,
\begin{equation}\label{eq:exactequiv}
F_p(x)\geq F_p(z)-\eps_p\quad(\forall x\in D)
\quad\Longleftrightarrow\quad
\de_{pj}(D;z)\leq\eps_{pj}\quad(j=1,\ldots,k_p).
\end{equation}

\end{theorem}
\begin{proof}
For each $j$,
\begin{align*}
&F_{pj}(x)\geq F_{pj}(z)-\eps_{pj}\quad(\forall x\in D)\\
&\iff \min_{x\in D}c_{pj}^{\mathsf T}x
\geq c_{pj}^{\mathsf T}z-\eps_{pj}\\
&\iff \de_{pj}(D;z)\leq\eps_{pj}.
\end{align*}

\end{proof}

The test requires one linear program for each protected criterion.
For a generic row $c\in\R^n$, write
$\de_c(D;z)=\max_{x\in D}c^{\mathsf T}(z-x)$. Using~\eqref{eq:Tineq},
the relevant program is
\begin{equation}\label{eq:primal}
\begin{aligned}
\max_x\quad &c^{\mathsf T}z-c^{\mathsf T}x\\
\text{subject to}\quad&\widehat Ax\leq\widehat b,\\
&\ell\leq x\leq u.
\end{aligned}
\end{equation}
Its dual also provides a direct certificate of the tolerance.

\begin{proposition}[Dual loss certificate]\label{prop:dual}
For $z\in D$, the exact loss equals
\begin{equation}\label{eq:dual}
\de_c(D;z)=c^{\mathsf T}z+
\min_{\substack{\lambda,\mu,\nu\geq0\\
\widehat A^{\mathsf T}\lambda+\mu-\nu=-c}}
\bigl(\widehat b^{\mathsf T}\lambda
+u^{\mathsf T}\mu-\ell^{\mathsf T}\nu\bigr).
\end{equation}
In particular, $\de_c(D;z)\leq\eps$ if and only if there exist
$\lambda,\mu,\nu\geq0$ satisfying
\begin{equation}\label{eq:dualfeas}
\begin{aligned}
\widehat A^{\mathsf T}\lambda+\mu-\nu&=-c,\\
c^{\mathsf T}z+\widehat b^{\mathsf T}\lambda
+u^{\mathsf T}\mu-\ell^{\mathsf T}\nu&\leq\eps.
\end{aligned}
\end{equation}
\end{proposition}
\begin{proof}
For any dual-feasible $\lambda,\mu,\nu$ and any $x\in D$,
\begin{align*}
c^{\mathsf T}(z-x)
&=c^{\mathsf T}z+\lambda^{\mathsf T}\widehat Ax
  +\mu^{\mathsf T}x-\nu^{\mathsf T}x\\
&\leq c^{\mathsf T}z+\widehat b^{\mathsf T}\lambda
  +u^{\mathsf T}\mu-\ell^{\mathsf T}\nu.
\end{align*}
The primal is feasible and bounded. Strong linear-programming duality
proves~\eqref{eq:dual}, and dual attainment gives~\eqref{eq:dualfeas}.
\end{proof}

There is also a sufficient test that requires no optimization. It ignores
the constraints of $T$ and uses only the coordinate bounds.

\begin{proposition}[Bound from the objective coefficients]\label{prop:coefficient}
If $z\in D$, then for every $c\in\R^n$,
\begin{equation}\label{eq:cheap}
0\leq\de_c(D;z)\leq
\sum_{i:c_i>0}c_i(z_i-\ell_i)
+\sum_{i:c_i<0}|c_i|(u_i-z_i).
\end{equation}
The right-hand side is the exact maximum over the full coordinate box.
If it does not exceed $\eps_{pj}$ for $c=c_{pj}$, the corresponding
criterion is protected on $D$.
\end{proposition}
\begin{proof}
For each $i$, the maximum of $c_i(z_i-x_i)$ over
$\ell_i\leq x_i\leq u_i$ equals
\[
\begin{cases}
c_i(z_i-\ell_i),&c_i>0,\\
|c_i|(u_i-z_i),&c_i<0,\\
0,&c_i=0.
\end{cases}
\]
Summing these maxima gives the bound. Restricting the full box to $D$
cannot increase it.
\end{proof}

For symmetric widths $|x_i-z_i|\leq r_i$, the simpler inequality is
\begin{equation}\label{eq:cheapsym}
\de_c(D;z)\leq\sum_{i=1}^n|c_i|r_i.
\end{equation}
This estimate may be conservative because it disregards the constraints
coupling the variables in $T$.

\begin{corollary}[Losses decrease under restriction]\label{cor:monotone}
If $D'\subseteq D$ are nonempty compact sets and $z\in\R^n$, then
$\de_c(D';z)\leq\de_c(D;z)$.
\end{corollary}
\begin{proof}
\[
\max_{x\in D'}c^{\mathsf T}(z-x)
\leq\max_{x\in D}c^{\mathsf T}(z-x).
\]
\end{proof}

\subsection{Choosing the width of a restriction}\label{sec:parametric}

The preceding results test given bounds. We now choose their width.
Let $D_0\subseteq T$ be a nonempty compact polytope containing $z$, with
$D_0=\{x\in\R^n:A_0x\leq b_0\}$. Initially $D_0=T$; at later levels
it also includes earlier restrictions. Fix positive scales
$s=(s_1,\ldots,s_n)\in\R_{+}^n$. For $r\geq0$, put
\begin{equation}\label{eq:Dr}
D(r)=\{x\in D_0:|x_i-z_i|\leq rs_i,\ i=1,\ldots,n\},
\end{equation}
and define the largest loss for a criterion row $c$ by
\begin{equation}\label{eq:deltafun}
\delta_c(r)=\max_{x\in D(r)}c^{\mathsf T}(z-x).
\end{equation}
Since $D(0)=\{z\}$ and $D(r_1)\subseteq D(r_2)$ for $r_1\leq r_2$,
these losses cannot decrease as the bounds widen.

\begin{theorem}[Properties of the maximum loss]
\label{thm:parametric}
For every $c\in\R^n$, the function
$\delta_c:[0,\infty)\to\R_+$ satisfies:

\begin{enumerate}
\item $\delta_c(0)=0$ and $\delta_c$ is increasing;
\item $\delta_c$ is continuous, concave and piecewise affine;
\item there exists $R<\infty$ such that
      $\delta_c(r)=\delta_c(R)$ for all $r\geq R$.
\end{enumerate}
\end{theorem}

\begin{proof}$\;$\\
\begin{enumerate}
\item By definition,
\[
D(0)=\{z\},\qquad
D(r_1)\subseteq D(r_2)
\quad\text{if }0\leq r_1\leq r_2.
\]
Consequently,
\[
\delta_c(0)=0,\qquad
\delta_c(r_1)\leq\delta_c(r_2).
\]

\item By linear programming duality,
\begin{equation}\label{eq:dualparametric}
\delta_c(r)=c^{\mathsf T}z+
\min_{(\lambda,\mu,\nu)\in Y_c}
\left[
b_0^{\mathsf T}\lambda
+z^{\mathsf T}(\mu-\nu)
+rs^{\mathsf T}(\mu+\nu)
\right],
\end{equation}
where
\[
Y_c=
\left\{
(\lambda,\mu,\nu)\geq0:
A_0^{\mathsf T}\lambda+\mu-\nu=-c
\right\}.
\]

The polyhedron $Y_c$ is nonempty and pointed.
Since the dual problem has a finite optimal value,
its minimum is attained at a vertex of $Y_c$.

Let $v^k=(\lambda^k,\mu^k,\nu^k)$,
$k=1,\ldots,N$, be the vertices of $Y_c$.
Then \eqref{eq:dualparametric} can be written as
\[
\delta_c(r)=\min_{1\leq k\leq N}
\{\alpha_k+\beta_k r\},
\]
where
\[
\begin{aligned}
\alpha_k
&=c^{\mathsf T}z+b_0^{\mathsf T}\lambda^k
  +z^{\mathsf T}(\mu^k-\nu^k),\\
\beta_k
&=s^{\mathsf T}(\mu^k+\nu^k).
\end{aligned}
\]

Thus, $\delta_c$ is the minimum of finitely
many affine functions. It follows that
$\delta_c$ is concave, continuous and piecewise affine.

\item Since $D_0$ is compact and $s_i>0$,
there exists $R<\infty$ such that
\[
|x_i-z_i|\leq Rs_i,
\qquad
\forall x\in D_0,\quad i=1,\ldots,n.
\]
Hence $D(r)=D_0$ for every $r\geq R$, and
\[
\delta_c(r)
=\max_{x\in D_0}c^{\mathsf T}(z-x)
=\delta_c(R).
\]
\end{enumerate}
\end{proof}
A smallest such saturation radius is
\begin{equation}\label{eq:R}
R=\max_{1\leq i\leq n}
\frac{\max\{\max_{x\in D_0}x_i-z_i,\ z_i-\min_{x\in D_0}x_i\}}{s_i}.
\end{equation}
It can be found by at most $2n$ linear programs.

Given protected criterion rows $c_1,\ldots,c_k$ and tolerances
$\eps\in\R_+^k$, the largest admissible width is
\begin{equation}\label{eq:rstar}
r^*(\eps)=\max\{r\in[0,R]:
\delta_{c_j}(r)\leq\eps_j,\ j=1,\ldots,k\}.
\end{equation}

\begin{proposition}[Largest admissible radius]
\label{prop:search}
Let $r^*(\eps)$ be defined by~\eqref{eq:rstar}.
Then the following properties hold:
\begin{enumerate}
\item The maximum $r^*(\eps)$ exists, and the set
of admissible radii is $[0,r^*(\eps)]$.

\item If $\eps\leq\eps'$, then
\[
r^*(\eps)\leq r^*(\eps').
\]

\item For any $0<\eta<R$, an admissible radius
$\widehat r$ satisfying
\[
0\leq r^*(\eps)-\widehat r\leq\eta
\]
can be obtained by bisection in at most
\begin{equation}\label{eq:count}
\left\lceil\log_2\frac{R}{\eta}\right\rceil
\end{equation}
iterations, each requiring at most $k$
linear programs.
\end{enumerate}
\end{proposition}

\begin{proof}$\;$\\
\begin{enumerate}
    \item 
 By Theorem~\ref{thm:parametric},
each $\delta_{c_j}$ is continuous and nondecreasing,
with $\delta_{c_j}(0)=0$.
Since $\eps_j\geq0$, the set
\[
\mathcal R(\eps)=
\left\{r\in[0,R]:
\delta_{c_j}(r)\leq\eps_j,\quad
j=1,\ldots,k\right\}
\]
is nonempty, closed and bounded.
Moreover, for $0\leq r_1\leq r_2$,
\[
r_2\in\mathcal R(\eps)
\quad\Longrightarrow\quad
r_1\in\mathcal R(\eps).
\]
Hence
\[
\mathcal R(\eps)=[0,r^*(\eps)],
\qquad
r^*(\eps)=\max\mathcal R(\eps).
\]

\item If $\eps\leq\eps'$, then
\[
\mathcal R(\eps)\subseteq\mathcal R(\eps').
\]
Consequently,
\[
r^*(\eps)\leq r^*(\eps').
\]
\item
 If $R$ is admissible, then
$r^*(\eps)=R$.
Otherwise, start with $a_0=0$ and $b_0=R$,
where $a_0$ is admissible and $b_0$ is not.

At each iteration, test the midpoint
\[
r_m=\frac{a_{m-1}+b_{m-1}}{2}
\]
and retain the half-interval containing $r^*(\eps)$.
Thus,
\[
a_m\leq r^*(\eps)\leq b_m,
\qquad
b_m-a_m=\frac{R}{2^m}.
\]
For
\[
m=\left\lceil\log_2\frac{R}{\eta}\right\rceil,
\]
we have
\[
0\leq r^*(\eps)-a_m
\leq\frac{R}{2^m}\leq\eta.
\]
Taking $\widehat r=a_m$ proves the result.
Each admissibility test requires at most
$k$ linear programs.
\end{enumerate}
\end{proof}

The one-dimensional search uses nested domains. If every coordinate
width were chosen independently, the resulting optimization problem
would not in general be convex. Bisection gives a certified accuracy,
not necessarily an exact breakpoint.

\subsection{Successive reductions and final guarantees}\label{sec:hierarchy}

We now apply the loss test at each decision level. Start with $D_1=T$.
For $p=1,\ldots,P-1$, choose a reference $z^{(p)}\in D_p$, tolerances
$\eps_p\in\R_+^{k_p}$, and bounds
$\ell^{(p)}\leq z^{(p)}\leq u^{(p)}$. Define
\begin{equation}\label{eq:recursive}
D_{p+1}=\{x\in D_p:\ell^{(p)}\leq x\leq u^{(p)}\}.
\end{equation}
Before passing $D_{p+1}$ to the next level, require
\begin{equation}\label{eq:stagecert}
\max_{x\in D_{p+1}}c_{pj}^{\mathsf T}(z^{(p)}-x)
\leq\eps_{pj},\qquad j=1,\ldots,k_p.
\end{equation}
The bounds can be chosen directly or by the width search of
Section~\ref{sec:parametric}, taking $D_0=D_p$.
A reference $z^{(p)}$ may be selected by maximizing
$w_p^{\mathsf T}F_p(x)$ on $D_p$ for $w_p>0$; the guarantees below
do not depend on this particular rule.

\begin{theorem}[Feasibility and retention at all earlier levels]\label{thm:hierarchy}
If the bounds contain $z^{(p)}$ and~\eqref{eq:stagecert} holds for
$p=1,\ldots,P-1$, then
\begin{equation}\label{eq:nested}
\varnothing\ne D_P\subseteq D_{P-1}\subseteq\cdots
\subseteq D_1=T\subseteq\calE.
\end{equation}
Every final choice $x^*\in D_P$ is efficient for each $F_p$ on $S$ and
satisfies
\begin{equation}\label{eq:finalguarantee}
F_p(x^*)\geq F_p(z^{(p)})-\eps_p,
\qquad p=1,\ldots,P-1.
\end{equation}
\end{theorem}
\begin{proof}
By construction, $z^{(p)}\in D_p$ and
$\ell^{(p)}\leq z^{(p)}\leq u^{(p)}$. Thus
$z^{(p)}\in D_{p+1}$ and $D_{p+1}\ne\varnothing$.
This proves~\eqref{eq:nested} by induction.
For $p<P$, $x^*\in D_P\subseteq D_{p+1}$, so~\eqref{eq:stagecert}
gives
\[
c_{pj}^{\mathsf T}(z^{(p)}-x^*)\leq\eps_{pj},
\qquad j=1,\ldots,k_p.
\]
Adding $d_{pj}$ proves~\eqref{eq:finalguarantee}.
Finally, $x^*\in T\subseteq\calE=\bigcap_{p=1}^P E_p$.
\end{proof}

\begin{corollary}[Earlier loss bounds remain valid]\label{cor:earlier}
For any $q\geq p+1$ and any criterion $j$ of level $p$,
\begin{equation}\label{eq:earlier}
\de_{pj}(D_q;z^{(p)})\leq
\de_{pj}(D_{p+1};z^{(p)})\leq\eps_{pj}.
\end{equation}
\end{corollary}
\begin{proof}
Use $D_q\subseteq D_{p+1}$ in Corollary~\ref{cor:monotone},
with the reference $z^{(p)}$ unchanged.
\end{proof}
A practical procedure is therefore: choose $z^{(p)}$ on $D_p$,
fix its tolerances, determine bounds passing~\eqref{eq:stagecert},
and pass $D_{p+1}$ to the next decision-maker. At the last level,
select $x^*\in D_P$. All earlier bounds remain in force.
The number of decision levels is finite; if widths are found by
bisection, the prescribed numerical accuracy determines the
number of loss tests.
\begin{remark}
Theorem \ref{thm:hierarchy} gives objective-retention guarantees on the precomputed
sorting polytope. It does not assert uniqueness or solve the nested
lower-level response problem.
\end{remark}

\subsection{Stability of the loss bounds}\label{sec:stability}

The exact loss changes when an objective coefficient or a restricted
set is perturbed. The following estimate covers both cases.
For nonempty compact sets $D,D'\subseteq\R^n$,
the Hausdorff distance induced by this norm is
defined by
\begin{equation}\label{eq:hausdorff}
d_H^\infty(D,D')
=
\max\left\{
\sup_{x\in D}\inf_{y\in D'}\|x-y\|_\infty,
\sup_{y\in D'}\inf_{x\in D}\|x-y\|_\infty
\right\}.
\end{equation}
We also use the dual norm
\[
\|c\|_1=\sum_{i=1}^{n}|c_i|,
\]
which satisfies
\[
|c^{\mathsf T}x|
\leq\|c\|_1\|x\|_\infty.
\]

\begin{proposition}[Stability of the worst loss]\label{prop:hausdorff}
Let $D,D'$ be nonempty compact sets, $z\in D$, $z'\in D'$, and
$c,c'\in\R^n$. Then
\begin{align}\label{eq:stab}
\bigl|\de_c(D;z)-\de_{c'}(D';z')\bigr|
\leq{}&\|c\|_1\bigl(\|z-z'\|_\infty
+d_H^\infty(D,D')\bigr)\notag\\
&+\|c-c'\|_1\max_{y\in D'}\|z'-y\|_\infty.
\end{align}
\end{proposition}

\begin{proof}
Set $h=d_H^\infty(D,D')$.

For every $x\in D$, there exists $y\in D'$
such that $\|x-y\|_\infty\leq h$.
By H\"older's inequality,
\begin{align*}
c^{\mathsf T}(z-x)
&=c^{\mathsf T}(z'-y)
  +c^{\mathsf T}(z-z')
  +c^{\mathsf T}(y-x)\\
&\leq \de_c(D';z')
  +\|c\|_1\|z-z'\|_\infty
  +\|c\|_1\|x-y\|_\infty\\
&\leq \de_c(D';z')
  +\|c\|_1(\|z-z'\|_\infty+h).
\end{align*}
Taking the maximum over $x\in D$ gives
\[
\de_c(D;z)-\de_c(D';z')
\leq \|c\|_1(\|z-z'\|_\infty+h).
\]
Similarly,
\[
\de_c(D';z')-\de_c(D;z)
\leq \|c\|_1(\|z-z'\|_\infty+h).
\]
Hence
\begin{equation}\label{eq:stab_sets}
|\de_c(D;z)-\de_c(D';z')|
\leq \|c\|_1(\|z-z'\|_\infty+h).
\end{equation}

Next, for every $y\in D'$,
\[
\begin{aligned}
&|c^{\mathsf T}(z'-y)
  -c'^{\mathsf T}(z'-y)|\\
&\qquad\leq
\|c-c'\|_1\|z'-y\|_\infty.
\end{aligned}
\]
Using the inequality
\[
\left|\max_{y\in D'}f(y)
-\max_{y\in D'}g(y)\right|
\leq\max_{y\in D'}|f(y)-g(y)|,
\]
we obtain
\begin{equation}\label{eq:stab_coeff}
|\de_c(D';z')-\de_{c'}(D';z')|
\leq\|c-c'\|_1R(D',z').
\end{equation}

Finally, the triangle inequality gives
\[
\begin{aligned}
&|\de_c(D;z)-\de_{c'}(D';z')|\\
&\leq|\de_c(D;z)-\de_c(D';z')|\\
&\quad+|\de_c(D';z')-\de_{c'}(D';z')|\\
&\leq\|c\|_1
  \bigl(\|z-z'\|_\infty+d_H^\infty(D,D')\bigr)\\
&\quad+\|c-c'\|_1R(D',z'),
\end{aligned}
\]
which proves~\eqref{eq:stab}.
\end{proof}

\begin{corollary}[Robust protection against coefficient errors]\label{cor:robust}
Suppose $\|\widetilde c-c\|_1\leq\rho$. If
\begin{equation}\label{eq:robust}
\de_c(D;z)+\rho\max_{x\in D}\|z-x\|_\infty\leq\eps,
\end{equation}
then
$\max_{x\in D}\widetilde c^{\mathsf T}(z-x)\leq\eps$
for every such $\widetilde c$.
\end{corollary}
\begin{proof}
Apply Proposition~\ref{prop:hausdorff} with $D'=D$ and $z'=z$.
\end{proof}

The Hausdorff estimate concerns the \emph{restricted} sets $D,D'$.
A small distance between the original polytopes alone need not give
the same estimate after both are intersected with fixed coordinate bounds.

\section{Numerical example}\label{sec:numerical}

We consider the two-level linear problem of \citet{KaciRadjef2023}:
\begin{equation}\label{eq:numproblem}
\left\{
\begin{aligned}
\text{Level }1:\quad
&\max_{x_1} F_1(x)=
\begin{pmatrix}
2x_1+2x_2\\
-\frac12 x_1+\frac7{25}x_2\\
-\frac15 x_1+\frac12 x_2
\end{pmatrix},\\
&\text{where }x_2\text{ solves}\\[2mm]
\text{Level }2:\quad
&\max_{x_2} F_2(x)=
\begin{pmatrix}
x_1+3x_2\\ -2x_1-x_2\\ x_2
\end{pmatrix},\\
&\text{subject to }x\in S.
\end{aligned}
\right.
\end{equation}
Here $S=\{x\in\R^2:Ax\leq b,\ x\geq0\}$, where
\[
A=\left(\begin{array}{rr}
-2&1\\ -1&2\\ 0&1\\ 1&0\\ -1&-2\\ 3&-4\\ 1&-2
\end{array}\right)\qquad
b=\left(\begin{array}{r}3\\9\\6\\6\\-9\\7\\2\end{array}\right).
\]
The common compromise set reported in \citet{KaciRadjef2023} is
\[
\mathcal E=
\operatorname{conv}\{(6,6),(3,6)\}\ \cup
\operatorname{conv}\{(3,6),(1,5)\}.
\]
We retain the sorting segment used in that paper,
\begin{equation}\label{eq:segment}
T=\operatorname{conv}\{(1,5),(3,6)\}
 =\{(1+2t,5+t):0\leq t\leq1\},
\end{equation}
which lies on $-x_1+2x_2=9$. The first-level reference is
$z=(2,5.5)$, corresponding to $t=1/2$. The adaptive method of
\citet{KaciRadjef2023} gives $x^K=(2.5,5.75)$ for its chosen
variable bounds. We do not repeat that method's construction of $T$;
we examine the objective losses caused by restrictions imposed on it.

Substitution into the objective functions gives
\begin{align}
F_1(t)&=(12+6t,\ 0.9-0.72t,\ 2.3+0.1t),\label{eq:F1t}\\
F_2(t)&=(16+5t,\ -7-5t,\ 5+t).\label{eq:F2t}
\end{align}

\paragraph{A first certified restriction.}
Take the tolerance vector
\begin{equation}\label{eq:epsnum}
\eps_1=(1.2,\ 0.144,\ 0.02).
\end{equation}
Since $F_1(z)=(15,0.54,2.35)$, the condition
$F_1(t)\geq F_1(z)-\eps_1$ is equivalent to
\[
12+6t\geq13.8,\qquad
0.9-0.72t\geq0.396,\qquad
2.3+0.1t\geq2.33.
\]
It follows that
\begin{equation}\label{eq:tnum}
D=\{(1+2t,5+t):0.3\leq t\leq0.7\},
\end{equation}
or, in terms of the decision variables,
\begin{equation}\label{eq:boundsnum}
1.6\leq x_1\leq2.4,\qquad 5.3\leq x_2\leq5.7.
\end{equation}
The exact worst losses on $D$ are
\begin{equation}\label{eq:lossesnum}
(\Delta_{11},\Delta_{12},\Delta_{13})=(1.2,\ 0.144,\ 0.02)=\eps_1.
\end{equation}
Each bound is attained at an endpoint of $D$.

For the second decision-maker, we maximize the equally weighted sum
$F_{21}(t)+F_{22}(t)+F_{23}(t)=14+t$ over $[0.3,0.7]$.
The selected decision is $t^*=0.7$, or $x^*=(2.4,5.7)$, with
\begin{equation}\label{eq:finalF}
F_1(x^*)=(16.2,\ 0.396,\ 2.37),\qquad
F_2(x^*)=(19.5,\ -10.5,\ 5.7).
\end{equation}
Its realized losses relative to $z$ are $(0,0.144,0)$ for $F_1$
and $(0,1,0)$ for $F_2$. In contrast, the original compromise
$x^K$ has a loss of $0.18$ in the second criterion of $F_1$.
Thus $x^K$ does not meet the tolerance $0.144$ chosen here;
the two methods use different restriction rules.

\paragraph{Twenty tolerance scenarios.}
To assess the effect of the tolerance level, we set
\begin{equation}\label{eq:gammatol}
\eps_1(\gamma)=\gamma(1.2,\ 0.144,\ 0.02),\qquad \gamma\geq0,
\end{equation}
and consider twenty values of $\gamma$ between $0$ and $3$.
For a common radius $0\leq r\leq1$, impose
$|x_1-2|\leq r$ and $|x_2-5.5|\leq r/2$. The retained set is
\begin{equation}\label{eq:tfromr}
D(r)=\left\{(1+2t,5+t):\frac{1-r}{2}\leq t\leq\frac{1+r}{2}\right\}.
\end{equation}
The worst first-level losses and the corresponding largest
admissible radius are
\begin{equation}\label{eq:exactradial}
\boldsymbol\Delta_1(D(r);z)=(3r,\ 0.36r,\ 0.05r),
\end{equation}
\begin{equation}\label{eq:radiuscompare}
r^*(\gamma)=\min\{0.4\gamma,1\}.
\end{equation}
Maximizing the second-level weighted sum then gives
\begin{equation}\label{eq:gammasol}
t^*(\gamma)=\frac{1+r^*(\gamma)}2,\qquad
x^*(\gamma)=\left(2+r^*(\gamma),\ 5.5+\frac{r^*(\gamma)}2\right).
\end{equation}

To separate the loss at the selected decision from the worst loss on
the entire retained set, write
\[
L_{pj}(\gamma)=\max\{0,F_{pj}(z)-F_{pj}(x^*(\gamma))\},\quad
\Delta_{pj}(\gamma)=\max_{x\in D(r^*(\gamma))}\bigl(F_{pj}(z)-F_{pj}(x)\bigr).
\]
With $r^*=r^*(\gamma)$, direct calculation yields
\begin{equation}\label{eq:realizedlosses}
L_1=(0,0.36r^*,0),\qquad L_2=(0,2.5r^*,0),
\end{equation}
whereas
\begin{equation}\label{eq:worstlosses}
\boldsymbol\Delta_1=(3r^*,0.36r^*,0.05r^*),\qquad
\boldsymbol\Delta_2=(2.5r^*,2.5r^*,0.5r^*).
\end{equation}
The prescribed bounds apply to $\boldsymbol\Delta_1$ only.
The second-level losses are reported relative to the same reference
$z$ for comparison.

The coefficient-only test uses the original matrix
\[
C_1=\begin{pmatrix}
2&2\\-1/2&7/25\\-1/5&1/2
\end{pmatrix}.
\]
It gives the bounds $(3r,0.64r,0.45r)$ and certifies at most
\begin{equation}\label{eq:roughradial}
r_{\rm coeff}(\gamma)=\min\left\{\frac{2\gamma}{45},1\right\}.
\end{equation}
At $\gamma=1$, this yields $r_{\rm coeff}=0.0444$, compared with
$r^*=0.4$. The exact certificate is less restrictive because it
takes the segment equation into account.

\paragraph{Results.}
Tables~\ref{tab:sensitivity} and~\ref{tab:losses} report the
selected decisions, maximum losses and measured computational cost.
The losses in Table~\ref{tab:losses} are worst-case values over the
whole reduced set, not losses realized at $x^*(\gamma)$.

\begin{table}[htbp]
\centering\normalsize
\caption{Tolerance scenarios, selected decisions and computational cost.
Times are medians of seven runs; RSS denotes peak process memory.}
\label{tab:sensitivity}
\setlength{\tabcolsep}{5pt}
\begin{tabular*}{\linewidth}{@{\extracolsep{\fill}}rrrrrrr@{}}
\toprule
$\gamma$&$r^*$&$r_{\rm coeff}$&$x_1^*$&$x_2^*$&Time (ms)&RSS (MiB)\\
\midrule
0.00 & 0.000 & 0.000 & 2.000 & 5.500 & 97.9 & 145.6 \\
0.10 & 0.040 & 0.004 & 2.040 & 5.520 & 104.3 & 145.7 \\
0.20 & 0.080 & 0.009 & 2.080 & 5.540 & 97.7 & 145.8 \\
0.30 & 0.120 & 0.013 & 2.120 & 5.560 & 105.5 & 145.8 \\
0.40 & 0.160 & 0.018 & 2.160 & 5.580 & 106.3 & 145.8 \\
0.50 & 0.200 & 0.022 & 2.200 & 5.600 & 104.6 & 145.8 \\
0.60 & 0.240 & 0.027 & 2.240 & 5.620 & 109.7 & 145.8 \\
0.70 & 0.280 & 0.031 & 2.280 & 5.640 & 107.1 & 145.8 \\
0.80 & 0.320 & 0.036 & 2.320 & 5.660 & 102.2 & 145.8 \\
0.90 & 0.360 & 0.040 & 2.360 & 5.680 & 100.9 & 145.8 \\
1.00 & 0.400 & 0.044 & 2.400 & 5.700 & 100.0 & 145.8 \\
1.25 & 0.500 & 0.056 & 2.500 & 5.750 & 102.2 & 145.8 \\
1.50 & 0.600 & 0.067 & 2.600 & 5.800 & 101.1 & 145.8 \\
1.75 & 0.700 & 0.078 & 2.700 & 5.850 & 100.4 & 145.8 \\
2.00 & 0.800 & 0.089 & 2.800 & 5.900 & 104.0 & 145.8 \\
2.25 & 0.900 & 0.100 & 2.900 & 5.950 & 101.7 & 145.8 \\
2.40 & 0.960 & 0.107 & 2.960 & 5.980 & 102.1 & 145.8 \\
2.50 & 1.000 & 0.111 & 3.000 & 6.000 & 3.9 & 145.8 \\
2.75 & 1.000 & 0.122 & 3.000 & 6.000 & 3.4 & 145.8 \\
3.00 & 1.000 & 0.133 & 3.000 & 6.000 & 3.6 & 145.8 \\
\bottomrule
\end{tabular*}
\end{table}
\begin{table}[htbp]
\centering\normalsize
\caption{Worst losses on the retained domain for both decision-makers.
Only the first-level losses are subject to prescribed tolerances.}
\label{tab:losses}
\setlength{\tabcolsep}{5pt}
\begin{tabular*}{\linewidth}{@{\extracolsep{\fill}}rrrrrrr@{}}
\toprule
$\gamma$&$\Delta_{11}$&$\Delta_{12}$&$\Delta_{13}$&
$\Delta_{21}$&$\Delta_{22}$&$\Delta_{23}$\\
\midrule
0.00 & 0.000 & 0.000 & 0.000 & 0.000 & 0.000 & 0.000 \\
0.10 & 0.120 & 0.014 & 0.002 & 0.100 & 0.100 & 0.020 \\
0.20 & 0.240 & 0.029 & 0.004 & 0.200 & 0.200 & 0.040 \\
0.30 & 0.360 & 0.043 & 0.006 & 0.300 & 0.300 & 0.060 \\
0.40 & 0.480 & 0.058 & 0.008 & 0.400 & 0.400 & 0.080 \\
0.50 & 0.600 & 0.072 & 0.010 & 0.500 & 0.500 & 0.100 \\
0.60 & 0.720 & 0.086 & 0.012 & 0.600 & 0.600 & 0.120 \\
0.70 & 0.840 & 0.101 & 0.014 & 0.700 & 0.700 & 0.140 \\
0.80 & 0.960 & 0.115 & 0.016 & 0.800 & 0.800 & 0.160 \\
0.90 & 1.080 & 0.130 & 0.018 & 0.900 & 0.900 & 0.180 \\
1.00 & 1.200 & 0.144 & 0.020 & 1.000 & 1.000 & 0.200 \\
1.25 & 1.500 & 0.180 & 0.025 & 1.250 & 1.250 & 0.250 \\
1.50 & 1.800 & 0.216 & 0.030 & 1.500 & 1.500 & 0.300 \\
1.75 & 2.100 & 0.252 & 0.035 & 1.750 & 1.750 & 0.350 \\
2.00 & 2.400 & 0.288 & 0.040 & 2.000 & 2.000 & 0.400 \\
2.25 & 2.700 & 0.324 & 0.045 & 2.250 & 2.250 & 0.450 \\
2.40 & 2.880 & 0.346 & 0.048 & 2.400 & 2.400 & 0.480 \\
2.50 & 3.000 & 0.360 & 0.050 & 2.500 & 2.500 & 0.500 \\
2.75 & 3.000 & 0.360 & 0.050 & 2.500 & 2.500 & 0.500 \\
3.00 & 3.000 & 0.360 & 0.050 & 2.500 & 2.500 & 0.500 \\
\bottomrule
\end{tabular*}
\end{table}

\begin{figure}[tbp]
\centering
\begin{minipage}[t]{0.49\linewidth}
\vspace{0pt}\centering
\includegraphics[width=\linewidth]{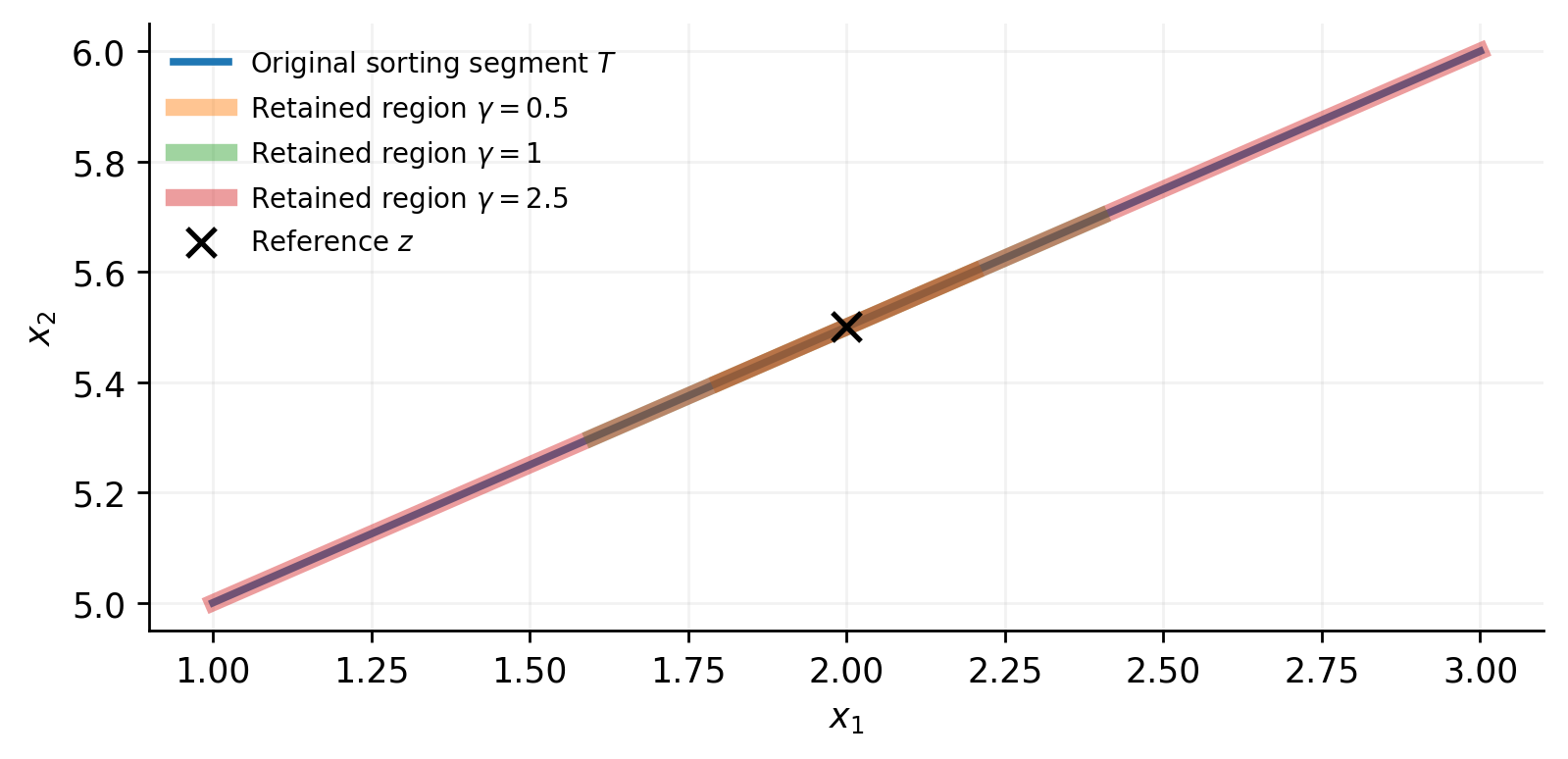}

\caption{Sorting segment and retained regions for selected tolerances.}
\label{fig:sorting_segment}
\end{minipage}
\hfill
\begin{minipage}[t]{0.49\linewidth}
\vspace{0pt}\centering
\includegraphics[width=\linewidth]{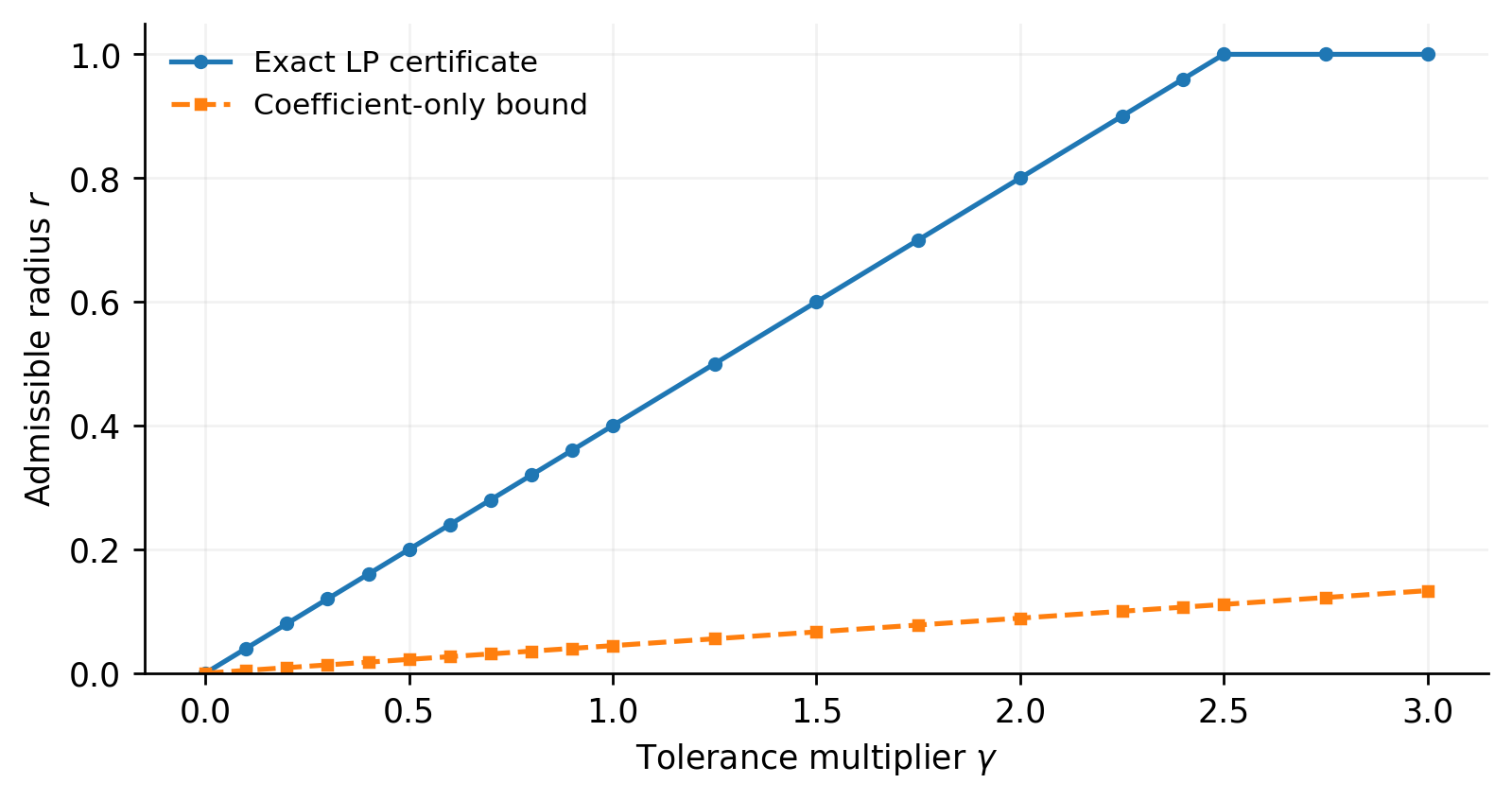}

\caption{Exact admissible radius and coefficient-only bound.}
\label{fig:radius}
\end{minipage}
\end{figure}

\begin{figure}[tbp]
\centering
\begin{minipage}[t]{0.49\linewidth}
\vspace{0pt}\centering
\includegraphics[width=\linewidth]{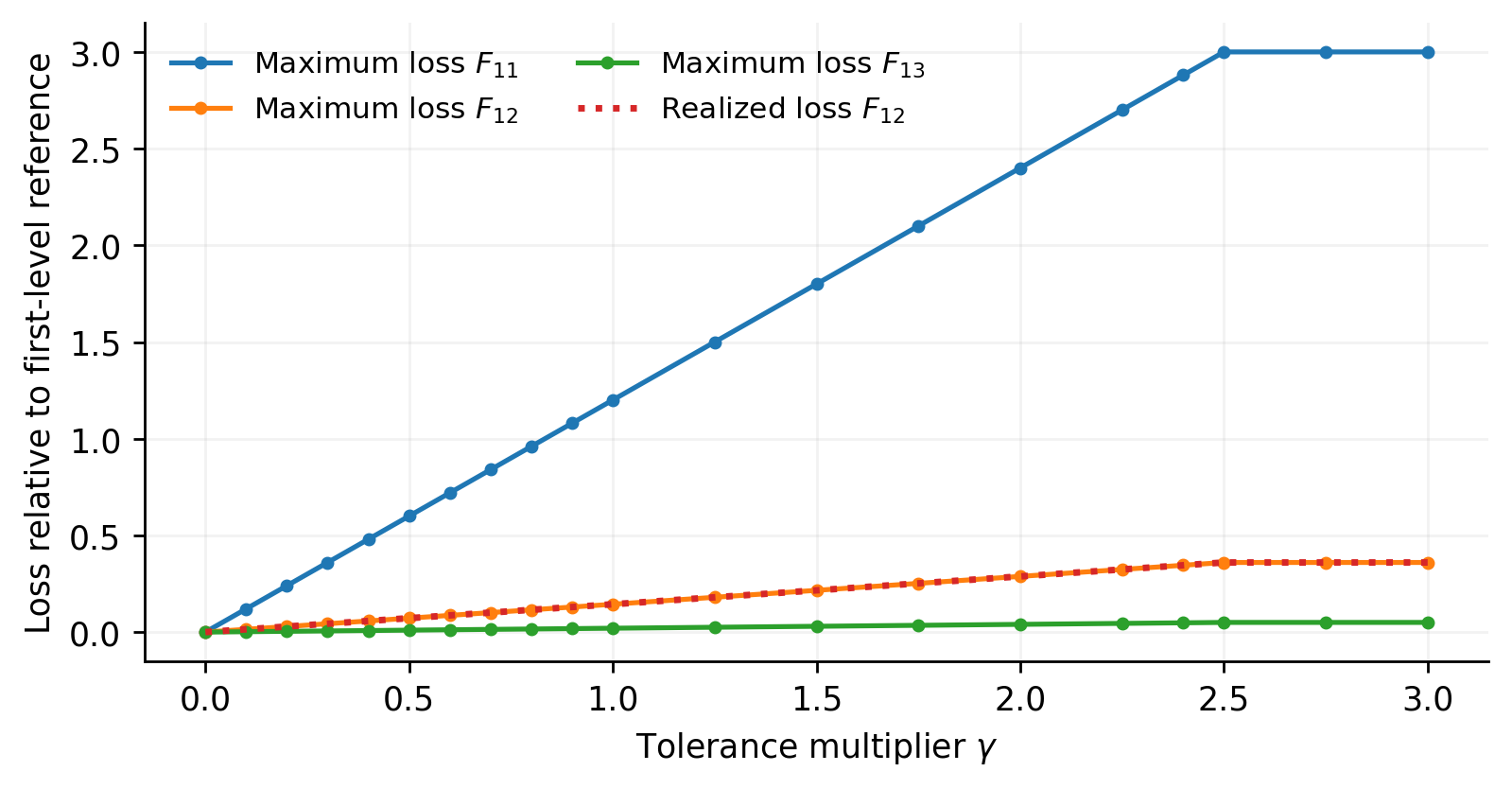}

\caption{Worst first-level losses and realized loss in the second criterion.}
\label{fig:loss_f1}
\end{minipage}
\hfill
\begin{minipage}[t]{0.49\linewidth}
\vspace{0pt}\centering
\includegraphics[width=\linewidth]{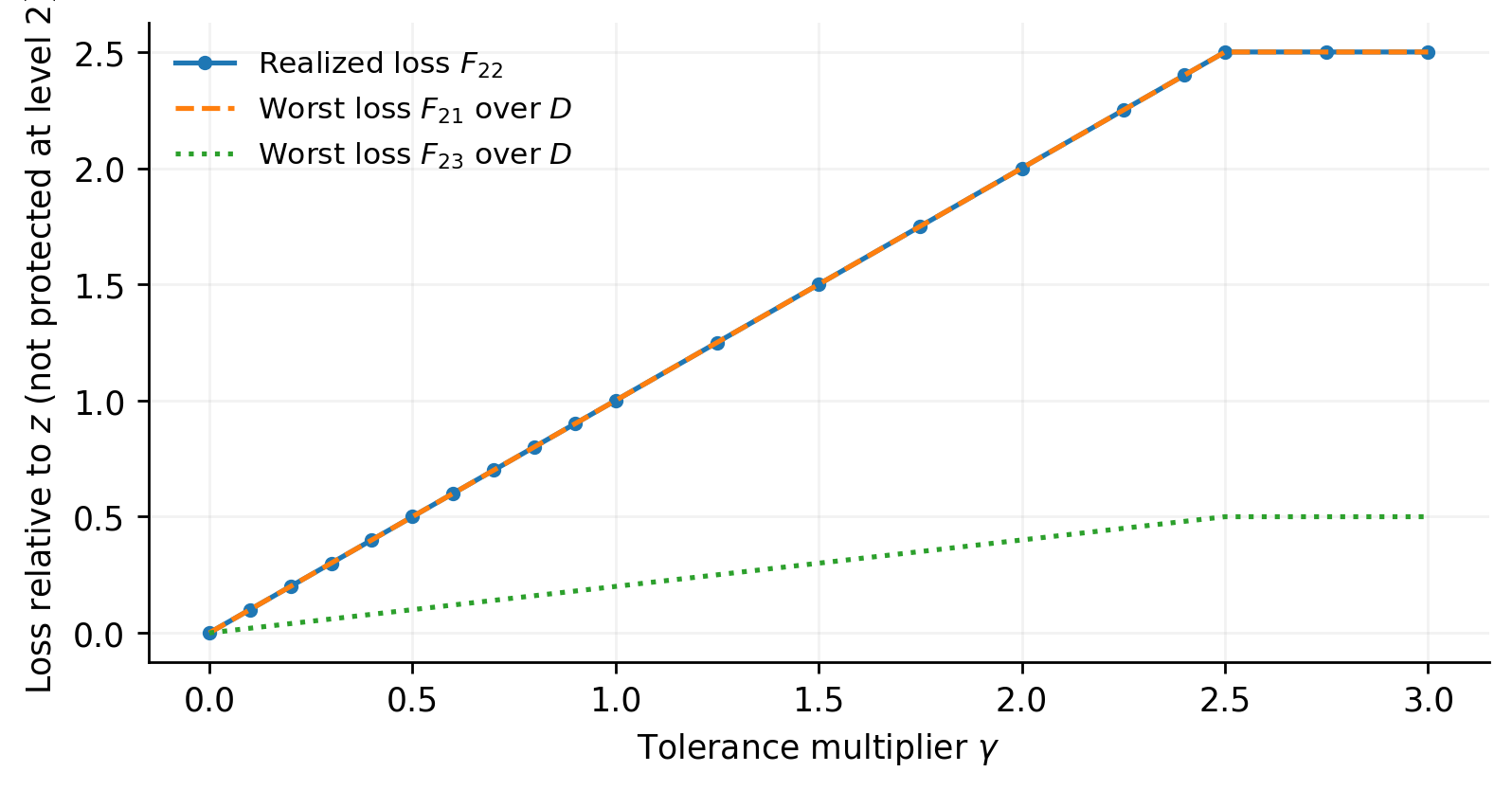}

\caption{Second-level losses relative to $z$, without a retention constraint.}
\label{fig:loss_f2}
\end{minipage}
\end{figure}

Figures~\ref{fig:sorting_segment}--\ref{fig:loss_f2} illustrate the
retained regions, admissible radii and objective losses.


\paragraph{Computational setting.}
We used Python~3.13.5, SciPy~1.17.0 (HiGHS) and
NumPy~2.3.5. Each radius test requires three linear programs.
The radius was found by bisection to an absolute precision of
$10^{-7}$. After warm-up, each scenario was repeated seven times.
We report the median elapsed time measured with
\texttt{time.perf\_counter} and the peak resident memory (RSS)
sampled with \texttt{psutil}~7.2.2. Plots were made with
Matplotlib~3.10.8.

The calculations were performed in a Linux container on a host
reporting an Intel Xeon Platinum 8573C processor. The container had
a quota of four CPU cores and 4~GiB of RAM; the full host
configuration was not exposed. Times cover the LP search, excluding
imports, plotting and construction of $T$. RSS refers to the whole
Python process, not only to the solver.

\begin{figure}[tbp]
\centering
\begin{minipage}[t]{0.49\linewidth}
\vspace{0pt}\centering
\includegraphics[width=\linewidth]{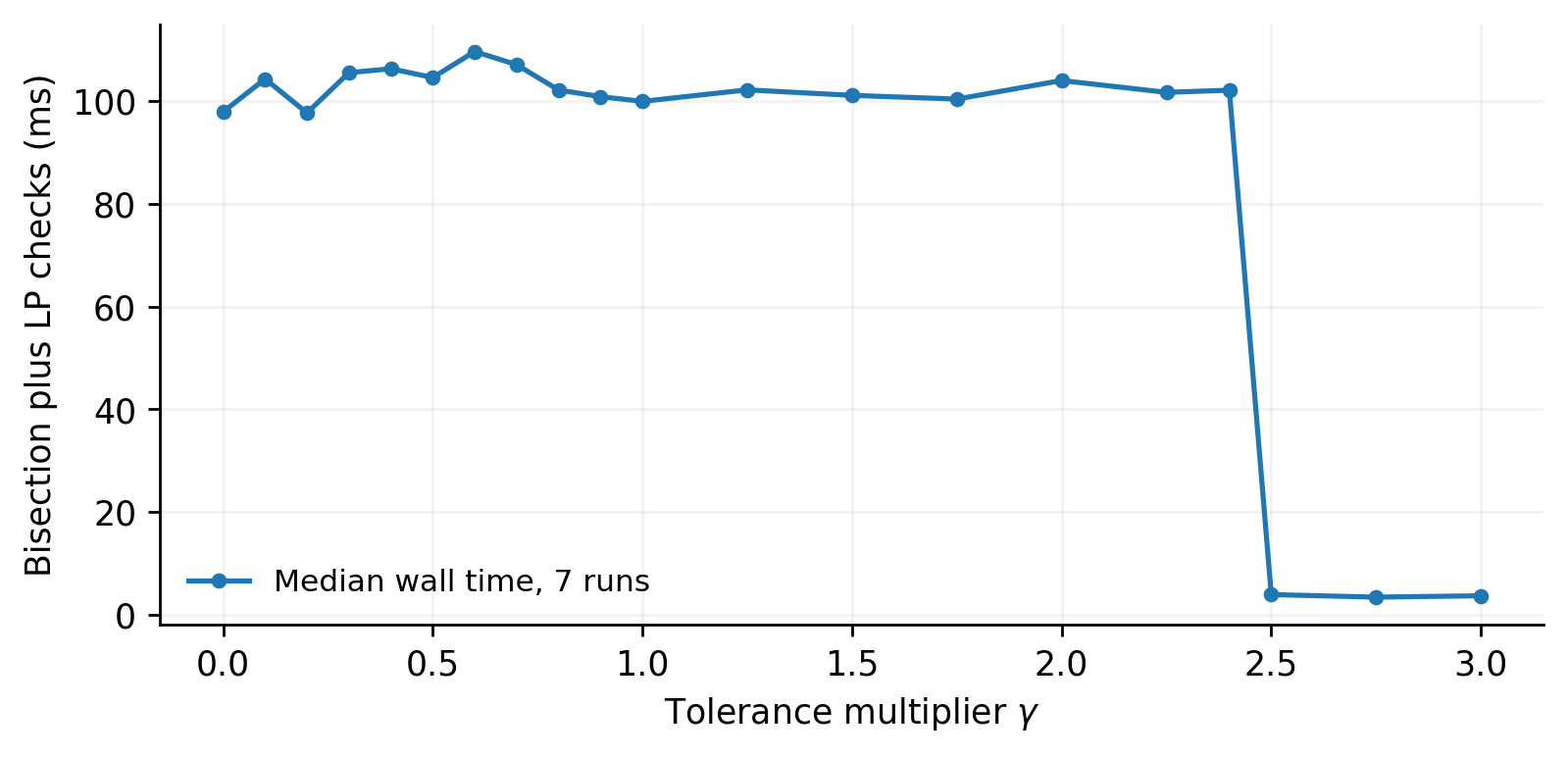}

\caption{Median solution time for each tolerance level (seven runs).}
\label{fig:time_profile}
\end{minipage}
\hfill
\begin{minipage}[t]{0.49\linewidth}
\vspace{0pt}\centering
\includegraphics[width=\linewidth]{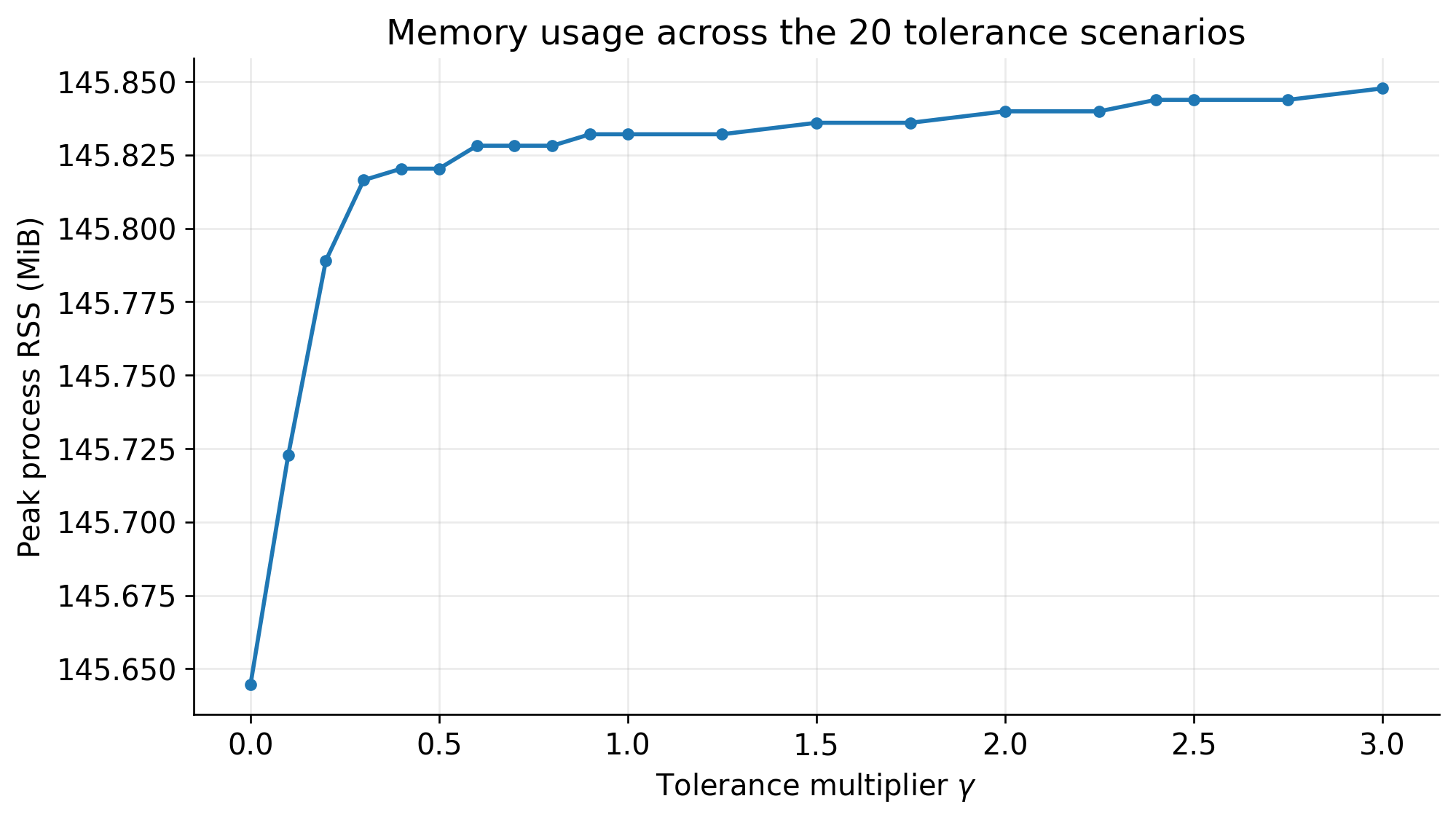}

\caption{Peak resident memory measured during the simulations.}
\label{fig:memory_profile}
\end{minipage}
\end{figure}
\begin{figure}[tbp]
\centering
\begin{minipage}[t]{0.49\linewidth}
\vspace{0pt}\centering
\includegraphics[width=\linewidth]{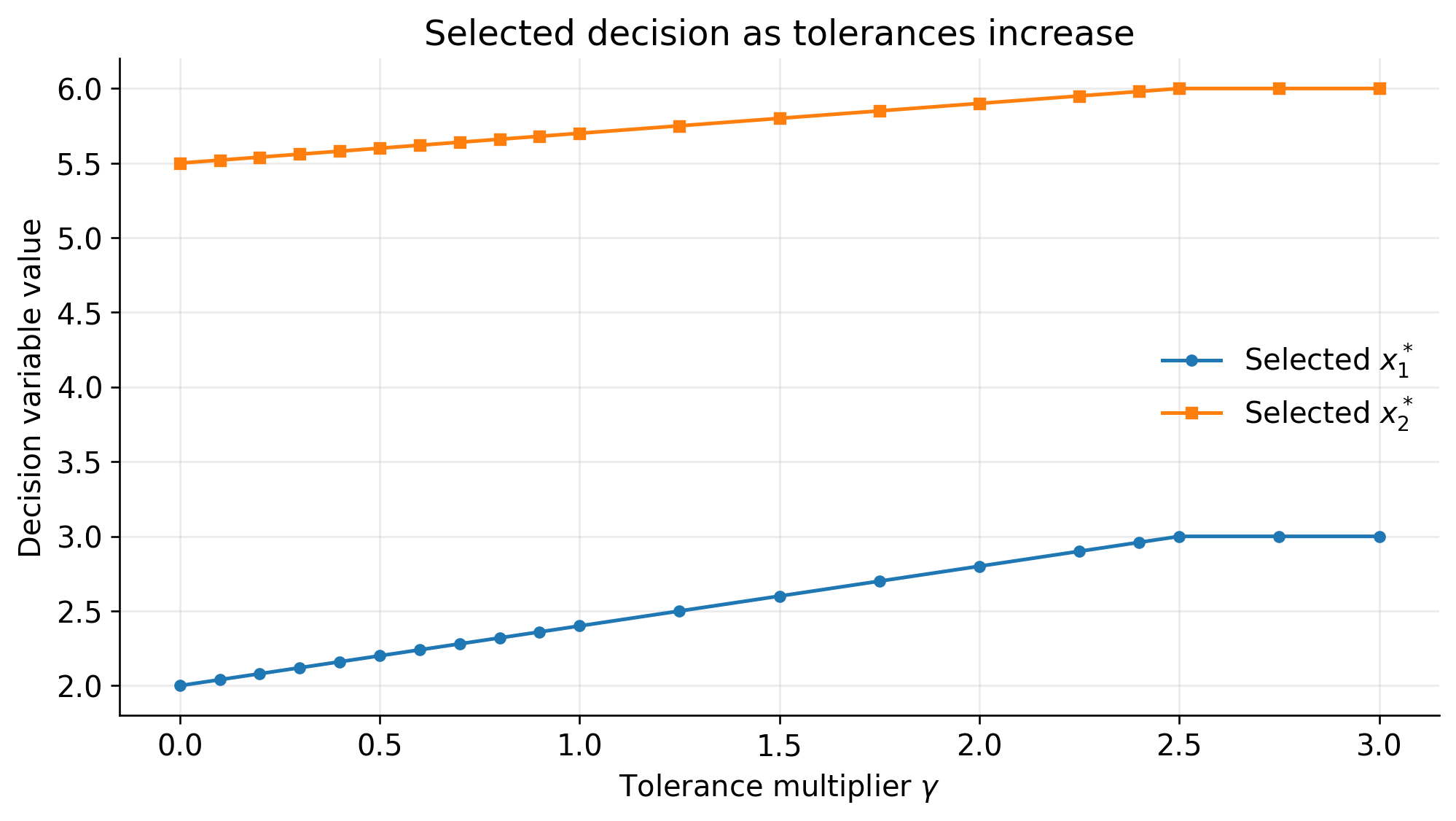}

\caption{Final decisions along the sorting segment as tolerances increase.}
\label{fig:selected_decisions}
\end{minipage}
\hfill
\begin{minipage}[t]{0.49\linewidth}
\vspace{0pt}\centering
\includegraphics[width=\linewidth]{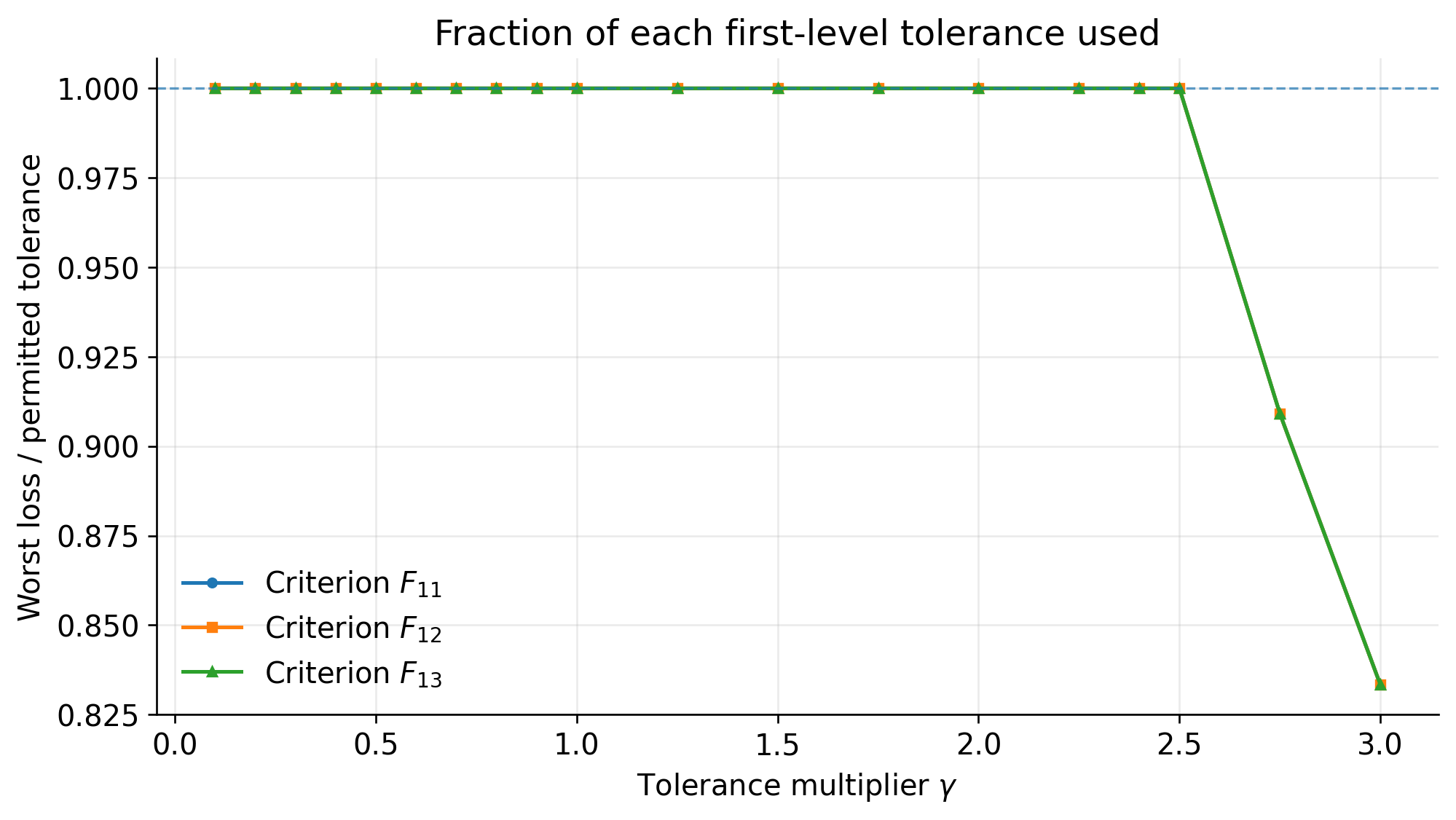}

\caption{First-level worst losses relative to their prescribed tolerances.}
\label{fig:tolerance_utilization}
\end{minipage}
\end{figure}

\paragraph{Discussion.}
The admissible radius increases linearly with $\gamma$ until
$\gamma=2.5$, when $D(r^*)=T$ and the selected point is $(3,6)$.
Beyond this value, larger tolerances do not change the result.
For $\gamma=1$, the maximum first-level loss vector is
$(1.2,0.144,0.02)$, while the realized loss at $(2.4,5.7)$ is
$(0,0.144,0)$. The second-level loss in its second criterion rises
with the selected radius; it is not subject to a retention bound.

For the seventeen scenarios with $\gamma<2.5$, the LP-based search
took roughly $100$--$110$~ms per run. At $\gamma\geq2.5$ the full
radius passes the first test, giving about $3$--$4$~ms. Peak process
memory was approximately $145.8$~MiB. These figures describe this
small problem and the stated computing environment; they do not
establish a general computational advantage.

\section{Conclusion}\label{sec:conclusion}

We have studied how variable bounds can preserve objective values during
successive reductions of a given sorting polytope. The exact loss is found
by linear programming, and its dual gives an algebraic certificate.
A simpler coefficient-based test is available, although it may exclude
restrictions that satisfy the exact tolerances. For proportional bounds,
the worst loss is continuous, nondecreasing, concave, and piecewise affine
in the width. This permits a certified search for the largest admissible
restriction.

When these bounds are imposed successively, the feasible sets remain
nonempty and every earlier tolerance remains valid for the final choice.
The stability estimate and two-level example make the distinction between
exact and coefficient-only guarantees explicit. The analysis presumes that
a sorting polytope has been constructed; numerical comparisons on larger
instances would be a useful next step.

\section*{Declaration of AI assistance}

AI-assisted tools were used to support the preparation
of this manuscript, including language editing,
LaTeX formatting, Python programming, and numerical
computations. The authors remain responsible for
the mathematical content, the verification of the
results, and the final version of the manuscript.

\end{document}